\documentclass[a4paper,oneside,10.5pt]{article}%
\usepackage{amsmath}
\usepackage{amsfonts}
\usepackage{amssymb}
\usepackage{graphicx}%
\providecommand{\U}[1]{\protect \rule{.1in}{.1in}}
\newtheorem{theorem}{Theorem}[section]

\newtheorem{corollary}[theorem]{Corollary}

\newtheorem{definition}[theorem]{Definition}
\newtheorem{assumption}[theorem]{Assumption}

\newtheorem{lemma}[theorem]{Lemma}

\newtheorem{remark}[theorem]{Remark}

\numberwithin{equation}{section}

\begin{document}

\title{The connection between discrete and continuous state constraints optimal control systems: the deterministic case}
\author{Shuzhen Yang\thanks{Zhong Tai Securities Institute for Financial Studies, Shandong University, PR China, (yangsz@sdu.edu.cn).}
\thanks{This work was supported by the National Natural Science Foundation of China (Grant No.11701330) and Young Scholars Program of Shandong University.}}
\date{}
\maketitle

\textbf{Abstract}: An optimal control problem driven by an ordinary differential equation under continuous state constraints is considered in this study. From an operational point of view, we introduce a discrete state constraints optimal control problem and prove that this discrete state constraints optimal control problem is a near-optimal control problem of the original problem. Furthermore, we show that the optimal solution of the near-optimal control problem converges to the optimal solution of the original one. Finally, we use a linear quadratic optimal problem to verify the main results of this study.


{\textbf{MSC2010}:  49K15, 34A60, 49N99.

\addcontentsline{toc}{section}{\hspace*{1.8em}Abstract}

\section{Introduction}
Since A.Ya. Dubovitskii and A.A. Milyutin  \cite{DM65}  proposed the optimality conditions for problems with state constraints, a vast amount of literature has been published on optimal control under state constraints in the field of optimality, and there are many related applications in economics and mathematical finance. For the deterministic optimal control problem case, we refer the reader to Frankowska \cite{F10} for a survey of the basic theory, such as optimal controls and the value function, and a discussion of the necessary optimality conditions under state constrained control systems. For different forms of maximum principle under state constraints, see the monograph by Vinter \cite{V00}, and for some historical comments on the maximum principle, see Pesch and Plail \cite{PP09}.

From the viewpoint of theory and applications of optimization, Zhou introduced the concept of near-optimization \cite{Z95,Z96}. Many authors have used the near-optimization method to study the switching linear quadratic (LQ) problem, stochastic recursive problem, linear forward backward stochastic systems, and so on; see \cite{LYZ05} for details.

In the deterministic optimal control problem, we use the following ordinal differential equation to describe the state,
\begin{equation}
\label{ine-1}
X(s)=x_0+\int_{0}^{s}b(X(t),u(t))dt,
\end{equation}
where $u(\cdot)$ is the control. In addition, we need to pay the running cost for the state and control, which we represent as $f(X(t),u(t))$ at time $t\in[0,T]$. In addition, we introduce the disposal cost for $X(\cdot)$ at time $T$, in general, and denote it by $\Psi(X(T))$. Thus, the cost functional is given as
\begin{equation}
\label{incos-1}
J(u(\cdot))=\displaystyle\int_0^Tf(X(t),u(t))dt+\Psi(X(T)).
\end{equation}
In reality, there are some limitations for state $X(\cdot)$ in $[0,T]$, i.e.,
\begin{equation}
\label{incons}
 X(t)\in \mathbb{S},\ \ 0\leq t\leq T,
\end{equation}
where $\mathbb{S}$ is a given set. For convenience, we call this problem the original one.

Frankowska \cite{F10} reviewed the classical maximum principle for the cost functional (\ref{incos-1}) under state constraints as follows: under some technical assumptions, the maximum principle under state constraints holds true. Consider the one-dimensional case and if $f=0$, there exist $\lambda\in \{0,1\}$, an absolutely continuous mapping
$p(\cdot) : [0,T] \to \mathbb{R}$, and a mapping $\phi(\cdot)$ of bounded total variation satisfying the adjoint equation
$$
dp(t)=b_x(\bar{X}(t),\bar{u}(t))(p(t)+\phi(t))
$$
and the maximum principle
$$
(p(t)+\phi(t))b(\bar{X}(t),\bar{u}(t))=\max_{u\in U}(p(t)+\phi(t))b(\bar{X}(t),u).
$$
and the transversality condition
$$
-p(T)-\phi(T)=\lambda \Phi(\bar{X}(T)).
$$
As Dmitruk pointed out in \cite{D09}, state constraints often appear in applied optimal control problems. However,
the solution of such problems is rather difficult because of a nonstandard form of
the adjoint equation in which almost nothing is known about the measure $\phi(\cdot)$.

In this study, our main idea is given as follows. From an operational point of view, we first introduce the following discrete state constraints,
\begin{equation}
\label{in0}
 X(t_i)\in \mathbb{S},\ \ i=0,1,\cdots,n,
\end{equation}
where $0=t_0< t_1<t_2<\cdots<t_n=T$. For the one-dimensional case, under some mild assumptions, we prove that the optimal control problem under discrete state constraints (\ref{in0}) is a near-optimal control problem of the original problem when $n$ is large enough.
In the end, we prove that the optimal solution of the optimal control problem under state constraints (\ref{in0}) converges to the optimal solution of the original problem.

In addition, there are many works related to Appendix A of this study. In deterministic case, A.V. Dmitruk and A.M. Kaganovich \cite{DK08,DK11} showed that this problem can be reduced to the standard problem of Pontryagin type by a simple change of variables, and
then one should only apply the classical maximum principle. In stochastic case, the stochastic maximum principle for a stochastic differential systems with a general cost functional (which is without state constraints) was developed in S. Yang \cite{Y16a},  the terminal cost functional is $\Psi(X_{[0,T]})$, where $X_{[0,T]}=X(s)_{0\leq s\leq T}$. However, there are some strong assumptions about Fr\'{e}chet
derivatives in \cite{Y16a}, and the structure of which is rather complicated; for further details see \cite{GY16,Y16b}. To remove some strong assumptions in \cite{Y16a}, the author investigated an optimal control problem with the following multi-time state cost functional,
\begin{equation}
J(u(\cdot))=E\big{[}\displaystyle\int_0^Tf(X(t),u(t))dt
+\Phi(X(\gamma_1),X(\gamma_2),\cdots,X(\gamma_N))\big{]},
\end{equation}
and a convex control domain $U$ in \cite{Y16c}.

The rest of this paper is organized as follows. In Section 2, we present the deterministic optimal control problem under state constraints. The connection between discrete and continuous state constraints optimal control problem is given in Section 4. In Section 4, we use an LQ problem to verify the main results of this study. In Appendix A, we develop the maximum principle for the cost functional (\ref{incos-1}) under multi-time state constraints (\ref{in0}).

\section{Optimal control problem}

Let $T>0$ be given, and consider the following controlled ordinary differential
equation,
\begin{equation}
\dot{X}^u(s)=b(X^u{(s)},u(s)),\quad s\in(0,T],\label{ODE_1}%
\end{equation}
with the initial condition $X^u(0)=x_0$, where
$u(\cdot)=\{u(s),s\in \lbrack0,T]\}$ is a control process taking value in a
 connected set $U$ of $\mathbb{R}^m$, which  generalize the convexity case,  and $b$ is a given deterministic function.

In this study, we consider the following cost functional,
\begin{equation}
J(u(\cdot))=%
{\displaystyle \int \limits_{0}^{T}}
f(X^u{(t)},u(t))dt+\Psi(X^u(T)),\label{cost-1}%
\end{equation}
under state constraints,
\begin{equation}
\label{cc-1}
 X^u(t)\in \mathbb{S},\ 0\leq t\leq T,
\end{equation}
where $\mathbb{S}$ is a given closed and connected set of $\mathbb{R}^m$ and
\[%
\begin{array}
[c]{l}%
b:\mathbb{R}^m\times U\to \mathbb{R}^m,\\
f:\mathbb{R}^m\times U\to \mathbb{R},\\
\Psi:\mathbb{R}^{m}\to \mathbb{R}.\\
\end{array}
\]
This is our original problem.

Let $b,f,\Psi$ be uniformly continuous and satisfy the following
 Lipschitz and continuous conditions.

\begin{assumption}
\label{ass-b}Suppose there exists a constant $c>0$ such that%
\[%
\begin{array}
[c]{c}%
\left| b(x_{1},u)-b(x_{2},u)\right|
 \leq c\left|x_1-x_2 \right|,\\
\end{array}
\]
$\forall(x_{1},u),(x_{2},u)\in{\mathbb{R}^m}\times U$.
\end{assumption}

\begin{assumption}
\label{assb-b2}
There exists a constant $c>0$ such that
$$
\displaystyle\sup_{u\in\mathcal{U}[0,T]}
\displaystyle\int_0^T\left|b(0,u(t))\right|^2dt \leq c,
$$
where $\mathcal{U}[0,T]=\{u(\cdot)\in L^2(0,T;U)\}.$
\end{assumption}

\begin{assumption}
\label{ass-fai}Let $b,f,\Psi$ be differentiable at $x$, and their derivatives in $x$ be continuous with respect to $(x,u)$.
\end{assumption}

\begin{assumption}
\label{ass-fai1}Suppose there exist constants $c_1,c_2>0$ such that%
\[%
\begin{array}
[c]{c}%
c_1\left|u_1-u_2 \right|\leq \left| b(x,u_1)-b(x,u_2)\right|\leq c_2\left|u_1-u_2 \right|
  ,\\
\end{array}
\]
$\forall(x,u_1),(x,u_2)\in{\mathbb{R}^m}\times U$.
\end{assumption}

\begin{remark}
Note that $b$ is continuous with respect to $(x,u)$ and Assumption \ref{ass-fai1} guarantees that  $b$ is a monotonic function at $u$ in the sense of each dimension of $u$.
\end{remark}

In addition, we make the following assumption about the control ability of system (\ref{ODE_1}) on the boundary of $\mathbb{S}$ that is used to prove the existence of near-optimal control for the original problem.

\begin{assumption}
\label{ass-fai2} For any given $0\leq t\leq s\leq T,u(\cdot)\in\mathcal{U}[0,T]$, $y\in \partial \mathbb{S}$, there exist two controls $u_1(\cdot),u_2(\cdot)\in \mathcal{U}[t,s]$ such that
$$
X^{u_1}(s)\leq y\leq X^{u_2}(s),
$$
where $X^{u_1}(t)=X^{u_2}(t)=y$, and where $\partial \mathbb{S}$ is the boundary set of $\mathbb{S}$.
\end{assumption}

Let Assumptions \ref{ass-b} and \ref{assb-b2} hold; then, there exists a unique
solution $X$ for equation (\ref{ODE_1}). Minimize (\ref{cost-1}) over $ \mathcal{U}[0,T]$ under constrained conditions (\ref{cc-1}); then, any $\bar{u}(\cdot)\in \mathcal{U}[0,T]$
satisfying
\begin{equation}
J(\bar{u}(\cdot))= \underset{u(\cdot)\in\mathcal{U}[0,T]}{\inf}J(u(\cdot)) \label{cost-2}%
\end{equation}
is called an optimal control. We denote this optimal problem as \textbf{OCSC}. The corresponding state trajectory $(\bar{u}(\cdot),\bar{X}(\cdot))$ is called an optimal state trajectory and optimal pair.

\section{The connection between discrete and continuous state constraints}
To deal with the problem \textbf{OCSC}, we will prove that the optimal control problem under the multi-time state constraints (\ref{cc-2}) provides a near-optimal control problem for the problem \textbf{OCSC}.
\subsection{Near-optimal control problem}
In this section, we consider the case $m=1$ and $\mathbb{S}=[0,+\infty)$ for the limitation of technique; for more cases see Remark \ref{cases2}. In the following, we first introduce the definition of near-optimal control for the optimal problem \textbf{OCSC}. For a given integer $n>0$ and $(t_0,t_1,\cdots,t_n)$ with $t_0=0,t_n=T$, let $(\tilde{u}^n(\cdot),\tilde{X}^n(\cdot))$ be an optimal pair of the cost functional (\ref{cost-1}) under the following state constraints,
\begin{equation}
\label{cc-2}
 X(t_i)\in \mathbb{S},\ i=0,1,\cdots,n.
\end{equation}

\begin{definition}
\label{den-1}
Suppose that $(\bar{u}^{}(\cdot),\bar{X}(\cdot))$ is an optimal pair of the problem \textbf{OCSC}. For any given $0<\delta<1$, if there exists an integer $n>0$ and $(t_0,t_1,\cdots,t_n)$, such that
$$
\left|J(\tilde{u}^{n}(\cdot))-J(\bar{u}(\cdot))\right|\leq \Theta(\delta),
$$
where $\Theta(\delta)\to 0$ as $\delta \to 0$. Then, we call that the cost functional (\ref{cost-1}) under state constraints (\ref{cc-2}) is a near-optimal control problem for the problem \textbf{OCSC}.
\end{definition}

Note that, under Assumptions \ref{ass-b} and \ref{assb-b2}, for any $u(\cdot)\in\mathcal{U}[0,T]$, the state equation (\ref{ODE_1}) admits a unique solution. We denote by
\begin{equation}
\label{ops-1}
\mathcal{A}[0,T]=\big{\{}(u(\cdot),X^u(\cdot)): X^u(\cdot) \text{ is the solution of equation (\ref{ODE_1})}, \text{ for } u(\cdot)\in \mathcal{U}[0,T]\big{\}},
\end{equation}
for $(u(\cdot),X^u(\cdot))\in \mathcal{A}[0,T]$, we define the norm on $\mathcal{A}[0,T]$ as
$$
\|(u(\cdot),X^u(\cdot))\|_{\mathcal{A}[0,T]}:=\displaystyle\sup_{0\leq t\leq T}\left|X^u(t)\right|+\displaystyle\int_0^T\left|u(t)\right|dt.
$$

In the following, we show some basic estimations for the state process that are related to the control $u(\cdot)\in\mathcal{U}[0,T]$.

\begin{lemma}
\label{lees-1}
Let Assumptions \ref{ass-b}, \ref{assb-b2} hold. For given $u(\cdot)\in\mathcal{U}[0,T]$, the solution $X^u(\cdot)$ of equation (\ref{ODE_1}) satisfies
\begin{equation}
\begin{array}
[c]{ll}
\displaystyle\sup_{0\leq t\leq T}\left|X^u(t)\right|\leq
L_1(x_0,T),\\
\displaystyle\sup_{t\leq r\leq s}\left|X^u(r)-X^u(t)\right|\leq L_2(x_0,T)(s-t),
\end{array}
\end{equation}
where $0\leq t\leq s\leq T$ and $L_1,L_2$ are the deterministic functions of $(x_0,T)$.
\end{lemma}
\noindent\textbf{Proof:} This lemma is classical; thus, we omit the proof. $\ \ \ \ \ \ \ \ \Box$

\bigskip

In the following, we give a lemma that describes the connected property of the set $\mathcal{A}[0,s]$ for $0\leq s\leq T$.

\begin{lemma}
\label{cu-1}
Suppose Assumptions \ref{ass-b} and \ref{assb-b2} hold. For any given $0\leq s\leq T$, $(u_1(\cdot),X^{u_1}(\cdot))$, \\
$ (u_2(\cdot),X^{u_2}(\cdot))$ $ \in \mathcal{A}[0,s]$, if $X^{u_1}(s)< X^{u_2}(s)$. Then, for any constant $K$ that satisfies $X^{u_1}(s)<K<X^{u_2}(s)$, there exists a control $u(\cdot)\in \mathcal{U}[0,s]$ such that $(u(\cdot),X^{u}(\cdot))\in \mathcal{A}[0,s]$ and
$$
X^u(s)=K.
$$
\end{lemma}
\noindent\textbf{Proof}: For given $0\leq s\leq T$, we suppose that the assertion of this lemma is incorrect. Then, there exists a constant $K_0$ satisfying $X^{u_1}(s)<K_0<X^{u_2}(s)$, and we have
\begin{equation}
\label{emp-1}
\mathcal{U}^0[0,s]=\big{\{}u(\cdot):X^u(s)=K_0   \big{\}}=\emptyset,
\end{equation}
where $\emptyset$ is an empty set. In the following, we denote by
$$
\mathcal{U}^{1}[0,s]=\big{\{}u(\cdot):X^u(s)\leq K_0 \big{\}}
$$
and
$$
\mathcal{U}^{2}[0,s]=\big{\{}u(\cdot):X^u(s)\geq K_0  \big{\}}.
$$
Obviously, $\mathcal{U}^{1}[0,s]$ and $\mathcal{U}^{2}[0,s]$ are not empty sets, indeed, $u_1(\cdot)\in \mathcal{U}^{1}[0,s]$ and $u_2(\cdot)\in \mathcal{U}^{2}[0,s]$.

By equality (\ref{emp-1}) and the definitions of $\mathcal{U}^{1}[0,s]$ and $\mathcal{U}^{2}[0,s]$, we can verify that they are closed sets and satisfy
$$
\mathcal{U}^{1}[0,s]\bigcap\mathcal{U}^{2}[0,s]=\emptyset
$$
and
$$
\mathcal{U}^{1}[0,s]\bigcup\mathcal{U}^{2}[0,s]=\mathcal{U}[0,s].
$$
For given $s$, note that $\mathcal{U}[0,s]$ is a connected set; thus, $\mathcal{U}^{1}[0,s]=\mathcal{U}^{2}[0,s]=\mathcal{U}[0,s]$, which is a contradiction.

This completes the proof. $\ \ \ \ \ \ \ \ \Box$

\bigskip

For given $0<\delta<1$, integer $n>0$, and $0=t_0< t_1< \cdots< t_n=T$, we denote by
$$
\mathcal{A}^0_n[0,T]=\big{\{}(u(\cdot),X^u(\cdot)): (u(\cdot),X^u(\cdot))\in\mathcal{A}[0,T], X^u(t_i)\in \mathbb{S},\ i=0,1,\cdots,n \big{\}}
$$
and
$$
\mathcal{A}^\delta_n[0,T]=\big{\{}(u(\cdot),X^u(\cdot)): (u(\cdot),X^u(\cdot))\in\mathcal{A}[0,T],X^u(t_i)\in \mathbb{S}_{\delta},\ i=0,1,\cdots,n \big{\}},
$$
where $\mathbb{S}_{\delta}=[-\delta,+\infty)$. Based on Lemma \ref{cu-1}, we have the following results.

\begin{lemma}
\label{cu-2}
Let Assumptions \ref{ass-b}, \ref{assb-b2}, \ref{ass-fai}, \ref{ass-fai1}, and \ref{ass-fai2} hold. If $\mathcal{A}^0_n[0,T]$ is not an empty set, then, for any $(u(\cdot),X^u(\cdot))\in \mathcal{A}^\delta_n[0,T]$ with $0<\delta<1$, there exists $(u^0(\cdot),X^{u^0}(\cdot))\in \mathcal{A}^0_n[0,T]$ and a constant $C>0$ such that
$$
\|(u(\cdot)-u^0(\cdot),X^u(\cdot)-X^{u^0}(\cdot))\|_{\mathcal{A}[0,T]}\leq C \sqrt{\delta}.
$$
\end{lemma}
\noindent\textbf{Proof.} Because $\delta>0$, by the definitions of $\mathcal{A}^0_n[0,T]$ and $\mathcal{A}^{\delta}_n[0,T]$, one obtains
$$
\mathcal{A}^0[0,T]\subset\mathcal{A}_n^{\delta}[0,T].
$$
If $(u(\cdot),X^u(\cdot))\in \mathcal{A}_n^0[0,T]$, we just take $(u^0(\cdot),X^{u^0}(\cdot))=(u(\cdot),X^u(\cdot))$. Thus, we only need to consider the case where
$$
(u(\cdot),X^u(\cdot))\in \mathcal{A}_n^{\delta,0}[0,T]=\mathcal{A}_n^{\delta}[0,T]-\mathcal{A}_n^0[0,T].
$$
In the following, we first prove that there exists $(u^0(\cdot),X^{u^0}(\cdot))\in \mathcal{A}_n^0[0,T]$ such that
 $$
\left| X^{u^0}(t_i)-X^{u}(t_i)\right|\leq \delta,\ i=0,1,\cdots,n.
 $$

\bigskip

 \textbf{Step 1}: Note that $\mathcal{A}_n^0[0,T]$ is not an empty set; thus, there exists $(u^1(\cdot),X^{u^1}(\cdot))\in \mathcal{A}_n^0[0,T]$, for $0\leq s\leq T$ such that
 $$
 X^{u^1}(s)=x_0+\displaystyle\int_0^{s}b(X^{u^1}(t),u^1(t))dt,
 $$
and $ X^{u^1}(t_i)\in \mathbb{S}$ with $0=t_0<t_1<\cdots<t_n=T$.

In the following, we construct the control $u^0(\cdot)$ step by step. We first consider time $t_1$; note that $(u(\cdot),X^{u}(\cdot))\in \mathcal{A}_n^{\delta,0}[0,T]$ satisfies
 $$
 X^{u}(t_1)=x_0+\displaystyle\int_0^{t_1}b(X^{u}(t),u(t))dt,
 $$
and $X^{u}(t_1)\in \mathbb{S}_{\delta}$. Note that $\mathbb{S}\subset \mathbb{S}_{\delta}$, if $X^{u}(t_1)\in \mathbb{S}$, we just take
$$
u^{11}(t)=u(t),\quad t\in[0,t_1);
$$
otherwise, $X^{u}(t_1)\in [-\delta,0)$, note that $X^{u^1}(t_1)\in \mathbb{S}$ and, by Lemma \ref{cu-1}, there exists $(u^{11}(\cdot),X^{u^{11}}(\cdot))\in \mathcal{A}_n^0[0,t_1]$ such that
$$
X^{u^{11}}(t_{1})=0,\ \left|X^{u^{11}}(t_{1})-X^{u}(t_1)\right|\leq \delta.
$$

Next, we consider time $t_2$. We consider the following four cases.

\textbf{Case 1:}
 $X^{u}(t_1),X^{u}(t_2)\in \mathbb{S}$, we just take
 $$
u^{22}(t)=u(t),\quad t\in[t_1,t_2).
$$
Thus,
$$
X^{u^{22}}(t_{2})=X^{u^{}}(t_{2}),\ \left|X^{u^{22}}(t_{2})-X^{u}(t_2)\right|\leq \delta.
$$

\textbf{Case 2:} $ X^{u}(t_1)\in \mathbb{S},X^{u}(t_2)\in [-\delta,0)$, because $X^{u}(\cdot)$ is continuous in $[t_1,t_2]$, there exists $s\in (t_1,t_2)$ such that
$$
X^{u}(s)=0.
$$
By Assumption \ref{ass-fai2}, there exist two controls $\tilde{u}^1(\cdot),\tilde{u}^2(\cdot)\in \mathcal{U}[s,t_2]$ such that
$$
X^{\tilde{u}^1}(t_2)\leq X^{u}(s)=0 \leq X^{\tilde{u}^2}(t_2).
$$
Again, by Lemma \ref{cu-1}, there exists a control $\tilde{u}^0(\cdot)\in \mathcal{U}[s,t_2)$ such that
$$
X^{\tilde{u}^0}(t_2)=0.
$$
Then, we denote by
$$
u^{22}(t)=u(t)1_{[t_1,s)}(t)+\tilde{u}^0(t)1_{[s,t_2)}(t),\ t\in [t_1,t_2),
$$
where
$$ 1_{[t_{1},s)}(t)=\left\{
\begin{aligned}
1,\quad t\in [t_{1},s)\\
0,\quad t\notin [t_{1},s)
\end{aligned}
\right.
\qquad
 1_{[s,t_2)}(t)=\left\{
\begin{aligned}
1,\quad t\in [s,t_2)\\
0,\quad t\notin [s,t_2).
\end{aligned}
\right.
$$
Therefore,
$$
X^{u^{22}}(t_{2})=0,\ \left|X^{u^{22}}(t_{2})-X^{u}(t_2)\right|\leq \delta.
$$

\textbf{Case 3:} $ X^{u}(t_1),X^{u}(t_2)\in [-\delta,0)$, similar to case 2, we can construct the control $u^{22}(\cdot)\in \mathcal{U}[t_1,t_2)$ such that
$$
X^{u^{22}}(t_{2})=0,\ \left|X^{u^{22}}(t_{2})-X^{u}(t_2)\right|\leq \delta.
$$

\textbf{Case 4:} $ X^{u}(t_1)\in [-\delta,0),X^{u}(t_2)\in \mathbb{S}$, by Lemma  \ref{cu-1}, there exists $u^{11}(\cdot)\in \mathcal{U}[t_0,t_1)$ such that $X^{u^{11}}(t_{1})=0$. We have
$$
X^{u}(t_2)=X^{u}(t_1)+\displaystyle\int_{t_1}^{t_2}b(X^u(s),u(s))ds
$$
and denote by $\hat{u}^0(t)=u^{11}(t)1_{[0,t_1)}(t)+u(t)1_{[t_1,t_2)}(t),\ t\in [0,t_2)$ and
$$
{X}^{\hat{u}^0}(t_2)=X^{u^{11}}(t_1)+\displaystyle\int_{t_1}^{t_2}b({X}^{\hat{u}^0}(s),u(s))ds.
$$
Because $X^{u}(t_1)<X^{u^{11}}(t_1)=0$, by a comparison with the theorem of ordinal differential equations, we obtain
\begin{equation}
\label{casee-1}
X^{u}(t_2)< {X}^{\hat{u}^0}(t_2).
\end{equation}
Again by Assumption \ref{ass-fai2} and an idea similar to that in case 2, we can construct a control $\hat{u}^1(\cdot)\in\mathcal{U}[0,t_2] $ that satisfies $\hat{u}^1(t)=u^{11}(t),\ t\in [0,t_1)$ such that
\begin{equation}
\label{casee-2}
X^{\hat{u}^1}(t_2)=0.
\end{equation}
Combining equations (\ref{casee-1}) and (\ref{casee-2}), and by Lemma \ref{cu-1}, we can obtain a control
$u^{22}(\cdot)\in \mathcal{U}[t_1,t_2)$ such that
$$
\left| X^{{u}^0}(t_i)-X^u(t_i) \right|<\delta, \quad   i=1,2,
$$
where
$$
u^0(t)=u^{11}(t)1_{[0,t_1)}(t)+u^{22}(t)1_{[t_1,t_2)}(t),\ t\in [0,t_2).
$$

Note that, when we consider time $t_3$, we only need to combine the value of $X^u(\cdot)$ at time $t_2$ and there are four cases similar to the time $t_2$. Thus, by mathematical induction, we can prove that there exists $(u^0(\cdot),X^{u^0}(\cdot))\in \mathcal{A}_n^0[0,T]$ and satisfying
 $$
 \left|X^{u^0}(t_i)-X^u(t_i)\right|\leq \delta,\ i=0,1,\cdots,n,
 $$
where $u^0(t)=\displaystyle\sum_{i=1}^nu^{ii}(t)1_{[t_{i-1},t_i)}(t)$.

 \bigskip

 \textbf{Step 2}: In this step, we prove that there exists $C>0$ such that
 $$
\|(u(\cdot)-u^0(\cdot),X^u(\cdot)-X^{u^0}(\cdot))\|_{\mathcal{A}[0,T]}\leq C \sqrt{\delta}.
$$
 Assumption \ref{ass-fai1} and Lemma \ref{lees-1} play a crucial role in this step. By Lemma \ref{lees-1}, we have
 $$
 \displaystyle\sup_{t\leq r\leq s}\left|X^u(r)-X^u(t)\right|\leq L_2(x,T)(s-t)
 $$
 and
 $$
 \displaystyle\sup_{t\leq r\leq s}\left|X^{u^0}(r)-X^{u^0}(t)\right|\leq L_2(x,T)(s-t).
 $$
In the following, let
 $$
 \Delta t_i=\frac{T}{n},
 $$
where $\Delta t_i=t_i-t_{i-1},\ i=1,2,\cdots,n$. Now, let $\displaystyle\frac{T}{n}=\sqrt{\delta}$ and note that $(u(\cdot),X^u(\cdot))\in \mathcal{A}_n^\delta[0,T]$ and $(u^0(\cdot),X^{u^0}(\cdot))\in \mathcal{A}_n^0[0,T]$, by Lemma \ref{lees-1}, we have
\begin{equation}
\label{eq-3}
\displaystyle\sup_{0\leq r\leq T}\left|X^{u^0}(r)-X^{u}(r)\right|\leq \displaystyle{C\sqrt{\delta}}.
\end{equation}
where $C>0$ is dependent on $x_0,T$, which will change line by line in the following. Recalling for any $0\leq s\leq T$,
 \begin{equation}
 \label{eq-1}
 X^{u^0}(s)=x_0+\displaystyle\int_0^{s}b(X^{u^0}(t),u^0(t))dt
 \end{equation}
and
 \begin{equation}
 \label{eq-2}
 X^{u}(s)=x_0+\displaystyle\int_0^{s}b(X^{u}(t),u(t))dt.
 \end{equation}
Combining equations (\ref{eq-1}) and (\ref{eq-2}), it follows that
\begin{equation}
\begin{array}
[c]{rl}
&X^{u^0}(t_i)-X^u(t_i)-X^{u^0}(t_{i-1})+X^u(t_{i-1})\\
=&\displaystyle\int_{t_{i-1}}^{t_i}\left(b(X^{u^0}(t),u^0(t))-b(X^{u}(t),u^0(t))\right)dt\\
&+\displaystyle\int_{t_{i-1}}^{t_i}\left(b(X^{u}(t),u^0(t))-b(X^{u}(t),u(t))\right)dt,
\end{array}
\end{equation}
For any given $1\leq i \leq n$,. By a simple calculation, we have
\begin{equation}
\label{inequa-11}
\begin{array}
[c]{rl}
\left|\displaystyle\int_{t_{i-1}}^{t_i}\left[b(X^{u}(r),u^0(r))-b(X^{u}(r),u(r))\right]dt\right|\leq C\delta.
\end{array}
\end{equation}
Note that $u(\cdot)\in \mathcal{U}[0,T]$ is measurable on $[0,T]$, for large $n$, we can choose a step function $u^n(\cdot)$ that is a constant on each interval $[t_{i-1},t_i),\ i=1,2,\cdots,n$, such that
\begin{equation}
\label{inequa-10}
\displaystyle\int_0^T\left| u^n(t)-u(t)\right|dt \leq C \delta.
\end{equation}
Paying attention to the details from Step 1, we can construct a control $u^0$ that is a constant on each interval $[t_{i-1},t_i)$ such that
$$
u^0(t)\leq u^n(t),\  \text{or}\ u^0(t)\geq u^n(t), \ t\in[t_{i-1},t_i),\ i=1,2,\cdots,n,
$$
from which we deduce that
$$
\displaystyle\int_{t_{i-1}}^{t_i}\left|b(X^{u}(r),u^0(r))-b(X^{u}(r),u^n(r))\right|dt
=\left|\displaystyle\int_{t_{i-1}}^{t_i}\left[b(X^{u}(r),u^0(r))-b(X^{u}(r),u^n(r))\right]dt\right|.
$$
Now, by Assumption \ref{ass-fai1} and inequalities (\ref{inequa-11}),(\ref{inequa-10}), we have
\begin{equation}
\label{inequa-12}
\begin{array}
[c]{rl}
& c_1\displaystyle\int_{t_{i-1}}^{t_i} \left| u^n(t)-u^0(t)\right|dt       \\
\leq& \left|\displaystyle\int_{t_{i-1}}^{t_i}\left[b(X^{u}(r),u^0(r))-b(X^{u}(r),u^n(r))\right]dt\right|\\
\leq &\left|\displaystyle\int_{t_{i-1}}^{t_i}\left[b(X^{u}(r),u^n(r))-b(X^{u}(r),u(r))\right]dt\right|+C\delta\\
\leq &C\delta.
\end{array}
\end{equation}
Thus, one obtains
$$
\left| u^n(t)-u^0(t)\right| \leq C \sqrt{\delta}, \quad t\in [0,T],
$$
it follows that
$$
\displaystyle\int_0^T\left| u(t)-u^0(t)\right|dt \leq C \sqrt{\delta},
$$
from which we can conclude that
$$
\|(u(\cdot)-u^0(\cdot),X^u(\cdot)-X^{u^0}(\cdot))\|_{\mathcal{A}[0,T]}\leq C \sqrt{{\delta}}.
$$
This completes the proof. $\ \ \ \ \ \ \ \ \Box$

\begin{theorem}
\label{th-cu}
Suppose that Assumptions \ref{ass-b}, \ref{assb-b2}, \ref{ass-fai}, \ref{ass-fai1}, and \ref{ass-fai2} hold. The problem \textbf{OCSC} exists as a near-optimal control problem.
\end{theorem}
\noindent\textbf{Proof}: Without loss of generality, we assume that the problem \textbf{OCSC} has an optimal pair $(\bar{u}(\cdot),\bar{X}(\cdot))$. In the following, we prove that there exist an integer $n>0$, $0=t_0<t_1< \cdots< t_n=T$, and an optimal pair $(\tilde{u}^n(\cdot),\tilde{X}^n(\cdot))\in \mathcal{A}_n^0[0,T]$ for the cost functional (\ref{cost-1}) under constrained conditions (\ref{cc-2}) such that
$$
\left|J(\tilde{u}^n(\cdot))-J(\bar{u}(\cdot))\right|\leq \Theta (\delta).
$$
 For given $0<\delta<1$, integer $n>0$, and $0=t_0< t_1< \cdots< t_n=T$, note that
$$
{\mathcal{A}}_n^{-\delta}[0,T]=\big{\{}(u(\cdot),X^u(\cdot)): (u(\cdot),X^u(\cdot))\in\mathcal{A}[0,T], \ X^u(t_i)\in \mathbb{S}_{-\delta},\ i=0,1,\cdots,n \big{\}},
$$
where $\mathbb{S}_{-\delta}=[\delta,+\infty)$  and
$$
{\mathcal{A}}_n^\delta[0,T]=\big{\{}(u(\cdot),X^u(\cdot)): (u(\cdot),X^u(\cdot))\in\mathcal{A}[0,T],\ X^u(t_i)\in \mathbb{S}_{\delta},\ i=0,1,\cdots,n \big{\}}
$$
and denote by
$$
\tilde{\mathcal{A}}[0,T]=\big{\{}(u(\cdot),X^u(\cdot)): (u(\cdot),X^u(\cdot))\in\mathcal{A}[0,T],\ X^u(t)\in \mathbb{S},\ 0\leq t\leq T \big{\}}.
$$
By Lemma \ref{lees-1}, we can choose a large enough $n$ such that
$$
{\mathcal{A}}_n^{-\delta}[0,T]\subset\tilde{\mathcal{A}}[0,T]\subset{\mathcal{A}}_n^\delta[0,T].
$$
Then, by Lemma \ref{cu-2}, for any $(u(\cdot),X^u(\cdot))\in {\mathcal{A}}_n^{\delta}[0,T]$, there exists $(u^n(\cdot),X^{u^n}(\cdot))\in {\mathcal{A}}_n^{-\delta}[0,T]$ such that
$$
\|(u(\cdot)-u^n(\cdot),X^u(\cdot)-X^{u^n}(\cdot))\|_{\mathcal{A}[0,T]}\leq C \sqrt{{\delta}},
$$
from which we deduce that for an optimal pair $(\tilde{u}^n(\cdot),\tilde{X}^n(\cdot))$ of the cost functional (\ref{cost-1}) under the following constrained conditions,
$$
\tilde{X}^{n}(t_{i})\in \mathbb{S},\ i=0,1,\cdots,n,
$$
there exists $(\tilde{u}(\cdot),\tilde{X}(\cdot))\in {\mathcal{A}}_n^{-\delta}[0,T]$ such that
\begin{equation}
\label{eq-5}
\|(\tilde{u}(\cdot)-\tilde{u}^n(\cdot),\tilde{X}(\cdot)-\tilde{X}^n(\cdot))\|_{\mathcal{A}[0,T]}\leq C \sqrt{{\delta}}.
\end{equation}
Recalling
$$
J(\tilde{u}(\cdot))=\big{[}\displaystyle\int_0^Tf(\tilde{X}(t),\tilde{u}(t))dt+\Psi(\tilde{X}(T))\big{]}
$$
and
$$
J(\tilde{u}^n(\cdot))=\big{[}\displaystyle\int_0^Tf(\tilde{X}^n(t),\tilde{u}^n(t))dt+\Psi(\tilde{X}^n(T))\big{]}.
$$
By equation (\ref{eq-5}) and Assumption \ref{ass-fai}, we have
\begin{equation}
\label{eqz-1}
0 \leq J(\tilde{u}(\cdot))-J(\tilde{u}^n(\cdot))\leq \Theta({\delta}).
\end{equation}

We assume that $(\bar{u}(\cdot),\bar{X}(\cdot))\in \tilde{\mathcal{A}}[0,T]$ is an optimal pair of the problem \textbf{OCSC}, from which we deduce that
\begin{equation}
\label{eqz-2}
J(\tilde{u}^n(\cdot)) \leq J(\bar{u}(\cdot))  \leq J(\tilde{u}(\cdot)).
\end{equation}
Then, combining equations (\ref{eqz-1}) and (\ref{eqz-2}), we obtain
$$
\left| J(\bar{u}(\cdot))-J(\tilde{u}^n(\cdot))  \right|\leq
\Theta(\delta),
$$
as $\delta\to 0$. This completes this proof.  $\ \ \ \ \ \ \ \ \Box$

\bigskip

Based on the results of Theorem \ref{th-cu}, there exist an integer $n>0$, $0=t_0< t_1< \cdots< t_n=T$, and an optimal pair $(\tilde{u}^n(\cdot),\tilde{X}^n(\cdot))$ for the cost functional (\ref{cost-1}) under constrained conditions (\ref{cc-2}) such that
$$
\left|J(\tilde{u}^n(\cdot))-J(\bar{u}(\cdot))\right|\leq \Theta(\delta).
$$

\begin{remark}
\label{cases2}
Recall from the proof of Lemma \ref{cu-2} that Assumption \ref{ass-fai1} is used to prove the approximation of the control. In fact, we can prove Theorem \ref{th-cu} without Assumption \ref{ass-fai1}, i.e., similar to the proof of Theorem \ref{th-cu}, the problem \textbf{OCSC} has a near-optimal control problem in the following two cases:

1) $f$ is independent of $u$;

2) $b$ and $f$ are linear functions of $u$.

\end{remark}

\subsection{Limitation of the near-optimal problem}
In the following, we show that the optimal solution of the near-optimal control problem converges to the optimal solution of the problem \textbf{OCSC}.

\begin{theorem}
\label{th-limit}
Suppose that Assumptions \ref{ass-b}, \ref{assb-b2}, \ref{ass-fai1}, and \ref{ass-fai2} hold. Let $f$ be independent of $u$ and $f,\Psi$ be strictly increasing on $\mathbb{R}$ with respect to $x$. There exists a sequence $\bar{X}^{\bar{u}^n}(\cdot)\in \mathbb{S}$ that converges to the optimal solution of the problem \textbf{OCSC}, where $(\bar{u}^{n}(\cdot),\bar{X}^{\bar{u}^n}(\cdot))$ is the optimal pair of the cost functional (\ref{cost-1}) under state constraints (\ref{cc-2}).
\end{theorem}
\noindent\textbf{Proof}: Suppose that $(\bar{u}^{n}(\cdot),\bar{X}^{\bar{u}^n}(\cdot))$ is the optimal pair of the cost functional as follows:
\begin{equation}
J(u(\cdot))=%
{\displaystyle \int \limits_{0}^{T}}
f(X^u{(t)})dt+\Psi(X^u(T))\label{cost-lim}%
\end{equation}
under state constraints
\begin{equation}
\label{lim-cst1}
0\leq X^u{(t_i)},\ i=0,1,\cdots,n,
\end{equation}
where $t_i-t_{i-1}=\frac{T}{n},\ i=1,2,\cdots,n$. Without loss of generality, we assume that $J(\bar{u}^n(\cdot))=0$ and set $\pi^n=\{ t_0,t_1,\cdots,t_n\}$. We add a new point $s_1=\frac{t_0+t_1}{2}$ in $\pi^n$ and denote the new division by $\pi^{n+1}=\{ s_0,s_1,\cdots,s_{n+1}\}$. Now, we consider the cost functional (\ref{cost-lim}) under the state constraints on $\pi^{n+1}$, i.e.,
\begin{equation}
\label{lim-cst2}
0\leq X^u{(s_i)},\ i=0,1,\cdots,n+1.
\end{equation}
Suppose that $(\bar{u}^{n+1}(\cdot),\bar{X}^{\bar{u}^{n+1}}(\cdot))$ is the optimal pair of the cost functional (\ref{cost-lim}) under the state constraints (\ref{lim-cst2}). It is easy to show that
$$
0=J(\bar{u}^n(\cdot))\leq  J(\bar{u}^{n+1}(\cdot)).
$$
We use the notation
$$
\tilde{X}^n(t)=\min(\bar{X}^{\bar{u}^n}(t),\bar{X}^{\bar{u}^{n+1}}(t)),\ 0\leq t\leq T.
$$
If there exists $h\in [0,T]$ such that $\tilde{X}^n(h)<\bar{X}^{\bar{u}^n}(h)$. Similar to the proof in Lemma \ref{cu-2}, we can construct a solution $\hat{X}^{\hat{u}}(\cdot)$ of equation (\ref{ODE_1}) such that
$$
\hat{X}^{\hat{u}}(t)\leq \bar{X}^{\bar{u}^n}(t),\ t\in[0,T]
$$
and $\hat{X}^{\hat{u}}(h)<\bar{X}^{\bar{u}^n}(h)$. Because $f$ and $\Psi$ are strictly increasing on $\mathbb{R}$ with respect to $x$, this is a contradiction. It follows that
$$
\bar{X}^{\bar{u}^n}(t)\leq \bar{X}^{\bar{u}^{n+1}}(t),\ t\in[0,T],
$$
by Assumption \ref{assb-b2}. Thus, the sequence $\{\bar{X}^{\bar{u}^n}(\cdot)\}_{n=1}^{\infty}$ is increasing and bounded and converges to the optimal solution of the problem \textbf{OCSC}.

This completes the proof.$\ \ \ \ \ \ \ \ \ \ \ \ \ \Box$

\begin{remark}
\label{cov-rek}
By Assumption \ref{assb-b2}, we obtain that the solution of state equation (\ref{ODE_1}) is bounded. We suppose that the solution of state equation (\ref{ODE_1}) takes value in $[x_{\min},x_{\max}]$. Thus, we can assume that $f,\Phi$ are strictly increasing on $[x_{\min},x_{\max}]$ in Theorem \ref{th-limit}.
\end{remark}

Based on Ekeland's variational principle, we have another limit result.

\begin{theorem}
\label{ek-2}
Let Assumptions \ref{ass-b}, \ref{assb-b2}, and \ref{ass-fai} hold. There exists a pair $(u^{n}(\cdot),X^{u^n}(\cdot))\in \mathcal{A}_n^{\frac{1}{n}}[0,T]$ that converges to the optimal pair of problem \textbf{OCSC} as $n\to \infty$.
\end{theorem}

\noindent\textbf{Proof.} Without loss of generality, we assume that $J(\bar{u}(\cdot))=0$, where $(\bar{u}(\cdot),\bar{X}(\cdot))$ is the optimal pair of problem (\ref{cost-3}) with constrained conditions (\ref{ccc-1}). For any integer $n>0$, we set
$$
J^{n}(u(\cdot))=\sqrt{\big{[}(J(u(\cdot))+\frac{1}{n})^+\big{]}^2
+\displaystyle\sum_{i=1}^n\big{[}(-X^u(t_i))^+\big{]}^2}.
$$
From Assumption \ref{ass-fai1}, one can verify that $J^{n}:\mathcal{U}[0,T]\to \mathbb{R}$ is continuous and satisfies
\begin{equation}
J^{n}(\bar{u}(\cdot))=\frac{1}{n}\leq \inf_{u\in\mathcal{U}[0,T]}J^{n}(u(\cdot))+\frac{1}{n}.
\end{equation}
Then, by Ekeland's variational principle, there exists a $u^{n}(\cdot)\in \mathcal{U}[0,T]$ such that
\begin{equation}
\label{ccci-evc}
J^{n}(u^{n}(\cdot))\leq J^{n}(\bar{u}(\cdot))=\frac{1}{n},\ \tilde{d}(u^{n}(\cdot),\bar{u}(\cdot))\leq \sqrt{\frac{1}{n}},
\end{equation}
where
$$
\tilde{d}(u^n(\cdot),\bar{u}(\cdot))=\sqrt{\int_0^T \left(u^n(t)-\bar{u}(t)\right)^2dt}.
$$
In addition, we have
\begin{equation*}
-\sqrt{\frac{1}{n}}\tilde{d}(u^{n}(\cdot),u(\cdot))\leq J^{n}(u(\cdot))-J^{n}(u^{n}(\cdot)),\ \forall u(\cdot) \in \mathcal{U}[0,T],
\end{equation*}
which deduces that
\begin{equation}
\label{cos-40c}
J^{n}(u^{n}(\cdot))+\sqrt{\frac{1}{n}}\tilde{d}(u^{n}(\cdot),u^{n}(\cdot))\leq J^{n}(u(\cdot))+\sqrt{\frac{1}{n}}\tilde{d}(u^{n}(\cdot),u(\cdot)),\ \forall u(\cdot) \in \mathcal{U}[0,T].
\end{equation}
Thus, inequality (\ref{cos-40c}) shows that $(u^{n}(\cdot),X^{u^n}(\cdot))$ is the optimal pair for the following cost functional
\begin{equation}
\label{cos-4c}
J^{n}(u(\cdot))+\sqrt{\frac{1}{n}}\tilde{d}(u^{n}(\cdot),u(\cdot)),
\end{equation}
without the state constraints.

By inequality (\ref{ccci-evc}), we have
$$
\sqrt{\big{[}(J(u^n(\cdot))+\frac{1}{n})^+\big{]}^2
+\displaystyle\sum_{i=1}^n\big{[}(-X^{u^n}(t_i))^+\big{]}^2}\leq \frac{1}{n},
$$
from which we can deduce that
$$
X^{u^n}(t_i)\geq -\frac{1}{n},\quad i=1,2,\ldots,n,
$$
and
$$
-\frac{2}{n}\leq J(u^n(\cdot))\leq 0.
$$
Thus, $(u^{n}(\cdot),X^{u^n}(\cdot))\in \mathcal{A}^{\frac{1}{n}}_n[0,T]$. This completes the proof.$\ \ \ \ \ \ \ \ \ \ \ \ \ \Box$

\section{LQ problem}
In this section, we consider the one-dimensional LQ optimal control problem that is used to verify the main results of this study.

Given the following linear state equation:
$$
X^u(t)=2+\int_0^t\left(X^u(s)+u(s)\right)ds,
$$
where $u(t)\in U=[-3,3],\ t\in [0,1]$ and the cost functional
\begin{equation}
\label{lq-cost}
J(u(\cdot))=\int_0^1\left(X^u(t)\right)^2dt+\left(X^u(1)\right)^2
\end{equation}
under state constraints
\begin{equation}
\label{lq-cons1}
1\leq X^u(t),\ t\in[0,1].
\end{equation}
Now, we use this model to verify Theorems \ref{ccc-th}, \ref{th-cu}, and \ref{th-limit}.

\bigskip

\textbf{Step 1:} In this step, we first verify Theorem \ref{ccc-th}. We introduce the following discrete version of the state constraints, and take $n=10$,
\begin{equation}
\label{lq-cons2}
1\leq X^u(t_i), \ i=0,1,\cdots,10.
\end{equation}
with $0=t_0<t_1<\cdots<t_{10}=1$ and $t_{i}-t_{i-1}=\frac{1}{10},\ i=1,2,\cdots,10$. We can verify that the optimal pair of the cost functional (\ref{lq-cost}) under state constraints (\ref{lq-cons2}) is given as follows,
 $$
(\bar{X}^{10}(t),\bar{u}^{10}(t))=\left\{
\begin{aligned}
(3-e^t,-3),\quad t\in [0,\ln2]\\
(3-2e^{t-\ln2},-3),\quad t\in (\ln2,l_7]\\
(-3+(6-2e^{l_7-\ln2})e^{t-l_7},3),\quad t\in (l_7,0.7]\\
(3-2e^{t-0.1*{(i-1)}},-3), \ i=8,9,10, \ t\in (0.1*(i-1),l_{i}],\\
 (-3+(6-2e^{l_i-0.1*{(i-1)}})e^{t-l_i},3),\ i=8,9,10, \ t\in (l_{i},0.1*i],
\end{aligned}
\right.
$$
 where $l_i,i=7,8,9,10$ is the unique solution of the following equation,
 $$
 -3+(6-2e^{l_7-\ln2})e^{0.7-l_7}=1
 $$
 and
  $$
-3+(6-2e^{l_i-0.1*{(i-1)}})e^{0.1*i-l_i}=1,\ i=8,9,10.
 $$

Let $p^{10}(\cdot)$ be the solution of the following series of first-order adjoint equations,
\begin{equation}%
\begin{array}
[c]{ll}%
&-d{p}^{10}(t)= \{p^{10}(t)-2\beta^0\bar{X}^{10}{(t)}\}dt,\ t\in(t_{i-1},t_{i}),\\
&p^{10}(t_{i})= -2\beta^0\bar{X}^{10}(t_{10})1_{i=10}(i)-\beta^i+p^{10}(t_{i}^{+}),  \quad i=1,2,\ldots,10,
\end{array}
\label{lqprin-1}%
\end{equation}
where $(\beta^0,\beta^1,\cdots,\beta^{10})$ satisfies
$$
\beta^0\geq 0,\ \ \left| \beta^0\right|^2+\displaystyle\sum_{j=1}^{10}\left| \beta^j\right|^2=1.
$$
Following the proof of Theorem \ref{ccc-th}, by a simple calculation, we obtain
$$
\beta^0=\beta^i=0,\beta^j\leq 0,\ i=7,8,9,10,\ j=1,2,3,4,5,6.
$$
and
$$
 p^{10}(t)=\left\{
\begin{aligned}
\left[\sum_{j=i}^6-\beta^je^{0.1*(j-i)}   \right]e^{t_i-t},\quad t\in (0.1*(i-1),0.1*i],\ i=1,2,3,4,5,6 \\
0,\quad t\in (0.6,1].
\end{aligned}
\right.
$$
The Hamilton function is as follows,
$$
H(\beta^0,x,u,p)=(x+u)p-\beta^0x^2.
$$
Thus,
$$
H(\beta^0,\bar{X}^{10}(t),\bar{u}^{10}(t),p^{10}(t))=(\bar{X}^{10}(t)+\bar{u}^{10}(t))p^{10}(t)-\beta^0\left(\bar{X}^{10}(t)\right)^2.
$$
Now, we plug the solution $(\bar{X}^{10}(\cdot),\bar{u}^{10}(\cdot))$ and $p^{10}(t)$ into the Hamilton function. It follows that
$$
H(\beta^0,\bar{X}^{10}(t),\bar{u}^{10}(t),p^{10}(t))\geq H(\beta^0,\bar{X}^{10}(t),u,p^{10}(t))
$$
which satisfies Theorem \ref{ccc-th}.

\bigskip

\textbf{Step 2:} In the following, we verify Theorem \ref{th-cu}. Note that the optimal pair of the cost functional (\ref{lq-cost}) under state constraints (\ref{lq-cons2}) is given as follows:
 $$
(\bar{X}(t),\bar{u}(t))=\left\{
\begin{aligned}
(3-e^t,-3),\quad t\in [0,\ln2]\\
(1,-1),\quad t\in (\ln2,1].
\end{aligned}
\right.
$$
In Step 1, we take $n=10$. It is easy to verify that
$$
\sup_{0\leq t\leq 1}\left|\bar{X}^n(t)-\bar{X}(t)\right|\leq \frac{C}{n},
$$
where $C$ is a constant that is independent of $n$. Thus, we obtain
$$
\left|J(\bar{u}^n(\cdot))-J(\bar{u}(\cdot))\right|< \frac{C}{n}.
$$

\bigskip

\textbf{Step 3:} By Remark \ref{cov-rek}, it is easy to verify the conditions of Theorem \ref{th-limit} for this LQ model. Following the proof of Theorem \ref{th-limit}, we add the new point in the division $0=t_0<t_1<\cdots<t_{10}=1$ step by step. By a simple calculation, we can obtain that
$$
\bar{X}^{10}(t)\leq \bar{X}^{11}(t),\ 0\leq t\leq 1,
$$
from which we can deduce that
$$
\bar{X}^n(t)\leq \bar{X}^{n+1}(t)\leq \bar{X}(t).
$$
By Step 2, we can obtain that $\sup_{0\leq t\leq 1}\left|\bar{X}^n(t)-\bar{X}(t)\right|$ converges to $0$ as $n\to \infty$.

\bigskip

\appendix

\section{The maximum principle for discrete state constraints}
In this section, we set $\mathbb{S}=[0,+\infty)$  and $m=1$ to simplify notation. We consider discrete constraints of (\ref{cc-1}), called multi-time state constraints, i.e.,
 \begin{equation}
\label{ccc-1}
0\leq  X^u(t_i),\ \ i=0,1,\cdots,n,
\end{equation}
 where $0=t_0<t_1<\cdots<t_n=T$. Next, we investigate the following general cost functional without constrained conditions
\begin{equation}
\label{costm-1}
\hat{J}(u(\cdot))=\displaystyle\int_0^Tf(X^u(t),u(t))dt+\psi(X^u(t_1),X^u(t_2),\cdots,X^u(t_n)),
\end{equation}
that is used to prove the maximum principle for cost functional (\ref{cost-2}) under the multi-time state constraints (\ref{ccc-1}). Note that we consider a general control domain $U$ that does not need to be convex. The main problem is to investigate the variational equation and adjoint equation. In the following, we introduce the first-order adjoint equations as follows,
\begin{equation}%
\begin{array}
[c]{rl}%
-d{p}(t)&= \{b_x(\bar{X}{(t)},\bar{u}(t))^{}p(t)
               -f_x(\bar{X}{(t)},\bar{u}(t))\}dt,\ t\in(t_{i-1},t_{i}),\\
p(t_{i}) &=-\psi_{x_i}(\bar{X}(t_1),\bar{X}(t_2),\cdots,\bar{X}(t_n))+p(t_{i}^{+}),\text{ \ }i=1,2,\cdots,n,
\end{array}
\label{prinm-1}%
\end{equation}
where $t_{i}^{+}$ is the right limit of $t_{i}$, and $p(t_{n}^{+})=0$.

We use the notation
\begin{equation*}
H(x,u,p)=b(x,u)^{}p-f(x,u),\text{ \  \ }%
\end{equation*}
where $(x,u,p)\in \mathbb{R}^m\times U\times \mathbb{R}^m $.

We have the following theorem.

\begin{theorem}
\label{Maximumprinciple} Let Assumptions \ref{ass-b}--\ref{ass-fai} hold,
and let $(\bar{u}(\cdot),\bar{X}(\cdot))$ be an optimal pair of (\ref{costm-1}).
Then there exists $p(\cdot)$ satisfying the series of first-order adjoint equations (\ref{prinm-1}) and such that
\begin{equation}%
\begin{array}
[c]{ll}%
&H(\bar{X}(t),\bar{u}(t),p(t))\geq H(\bar{X}(t),u,p(t))),
\end{array}
\label{prinm-3}%
\end{equation}
for any $u\in U$ and $t \in(t_{i-1},t_i)$, $i=1,2,\cdots,n$.
\end{theorem}

The proof of Theorem \ref{Maximumprinciple} is similar to that in the classical case; therefore, we could use the so-called spike
variation technique. Let $(\bar{u}(\cdot),\bar{X}(\cdot))$ be the given
optimal pair of cost functional (\ref{costm-1}). Let $\varepsilon>0,$ and $E_{\varepsilon}=[v,v+\varepsilon
]\subset (t_{i-1},t_i)$ for some $i\in\{1,2,\cdots,n\}$. Let $u(\cdot)\in \mathcal{U}[0,T]$ be any given control. We define the following%
\[
u^{\varepsilon}(t)=\left \{
\begin{array}
[c]{cl}%
\bar{u}(t), & \text{if }t\in \lbrack0,T]\backslash E_{\varepsilon},\\
u(t), & \text{if }t\in E_{\varepsilon},
\end{array}
\right.
\]
where, obviously, $u^{\varepsilon}(\cdot)\in \mathcal{U}[0,T]$.  The following lemma is
useful for proving Theorem \ref{Maximumprinciple}.

\begin{lemma}
\label{lem-2} Let Assumptions \ref{ass-b} and \ref{assb-b2} hold, let
$X^{\varepsilon}(\cdot)$ be the solution of equation (\ref{ODE_1}) under the
control $u^{\varepsilon}(\cdot)$, and let $y^{}(\cdot)$ be the solution
of the following equations:%
\begin{equation}
\label{aprom-1}
\begin{array}
[c]{cl}%
d{y}(t)= &\big{[} b_x(\bar{X}{(t)},\bar{u}(t))y(t)
+b(\bar{X}{(t)},u^{\varepsilon}(t))-b(\bar{X}{(t)},\bar{u}(t))\big{]}dt, \\
y(0)= & 0,\quad t\in (0,T].
\end{array}
\end{equation}
Then%
\begin{equation}%
\begin{array}
[l]{l}%
\max_{t\in \lbrack0,T]}\left \vert y^{}(t)\right \vert
=O(\varepsilon^{}),\\
\max_{t\in \lbrack0,T]}\left \vert X^{\varepsilon}(t)-\bar{X}(t)-y^{
}(t)\right \vert =o(\varepsilon),  \\
 \end{array}
 \label{varm-1}%
\end{equation}
and%
\begin{equation}%
\begin{array}
[c]{rl}
& \hat{J}(u^{\varepsilon}(\cdot))-\hat{J}(\bar{u}(\cdot))\\
= &
\displaystyle \sum \limits_{i=1}^{n}\psi_{x_i}(\bar{X}(t_1),\bar{X}(t_2),\cdots,\bar{X}(t_n))y(t_{i}) \\
& +
E{\displaystyle \int \limits_{0}^{T}}
\{f_x(\bar{X}{(t)},\bar{u}(t))y(t)+f(\bar{X}{(t)},u^{\varepsilon}(t))-f(\bar{X}{(t)},\bar{u}(t))\}dt+o(\varepsilon).
\end{array}
\label{varm-3}%
\end{equation}
\end{lemma}

\medskip

\noindent
\textbf{Proof:}
By a simple calculation, we can prove equation (\ref{varm-1}). Note that%
\begin{equation}%
\begin{array}
[c]{rl}
& \hat{J}(u^{\varepsilon}(\cdot))-\hat{J}(\bar{u}(\cdot))\\
= & \psi(X^{\varepsilon}(t_1),X^{\varepsilon}(t_2),\cdots,X^{\varepsilon}(t_n))
-\psi(\bar{X}(t_1),\bar{X}(t_2),\cdots,\bar{X}(t_n))\\
&+{\displaystyle \int \limits_{0}^{T}}
\{f(X{}^{\varepsilon}(t),u^{\varepsilon}(t))-f(\bar{X}{(t)},\bar{u}(t))\}dt.\\
\end{array}
\label{valuem-0}%
\end{equation}
By equation (\ref{varm-1}), it follows that
\begin{equation}%
\begin{array}
[c]{rl}
& \hat{J}(u^{\varepsilon}(t))-\hat{J}(\bar{u}(t))\\
= & \displaystyle\sum_{i=1}^{n}\psi_{x_i}(\bar{X}(t_1),\bar{X}(t_2),\cdots,\bar{X}(t_n))y(t_i)\\
&+{\displaystyle \int \limits_{0}^{T}}
\{f_x(\bar{X}{(t)},\bar{u}(t))y(t)+f(\bar{X}{(t)},u^{\varepsilon}(t))-f(\bar{X}{(t)},\bar{u}(t))\}dt+o(\varepsilon).\\
\end{array}
\label{valuem-1}%
\end{equation}
This completes the proof. $\ \ \ \ \ \ \ \ \Box$

\bigskip

Similar to the techniques used in the proof of the deterministic maximum principle on each time interval $(t_{i-1},t_i)$ for $i=1,2,\cdots,n$. We now prove Theorem \ref{Maximumprinciple}.

\textbf{Proof of Theorem \ref{Maximumprinciple}}. For $\ t\in (
t_{i-1},t_{i}),$ applying the differential chain rule to $p(t)^{}y^{
}(t)$, and by Assumption \ref{ass-fai}, we have%
\begin{equation}%
\begin{array}
[c]{rl}
& p(t_{i})^{}y(t_{i})-p(t_{i-1}^{+})^{}y^{}(t_{i-1})\\
= & -\psi_{x_i}(\bar{X}(t_1),\bar{X}(t_2),\cdots,\bar{X}(t_n))y(t_{i})+p(t_{i}^{+})y^{}(t_{i})
-p(t_{i-1}^{+})y^{}(t_{i-1}^{})\\
= &\displaystyle \int \limits_{t_{i-1}}^{t_{i}}
\big{[}f_x(\bar{X}{(t)},\bar{u}(t))y_{}
^{}(t)+p(t)^{}(b(\bar{X}{(t)},u^{\varepsilon}(t))-b(\bar{X}{(t)},\bar{u}(t)))\big{]}dt.
\end{array}
\label{maxm-1}%
\end{equation}
Adding $i$ to both sides of equation (\ref{maxm-1}), it follows that
\[%
\begin{array}
[c]{rl}
&
\displaystyle \sum \limits_{i=1}^{n}
\big{[} p(t_{i})y(t_{i})-p(t_{i-1}^{+})y^{}(t_{i-1})\big{]} \\
= &
{\displaystyle \sum \limits_{i=1}^{n}}
\{-\psi_{x_i}(\bar{X}(t_1),\bar{X}(t_2),\cdots,\bar{X}(t_n))y(t_{i})+p(t_{i}^{+})y^{}(t_{i})
-p(t_{i-1}^{+})y^{}(t_{i-1}^{}))\} \\
= &  \displaystyle \sum \limits_{i=1}^{n}-\psi_{x_i}(\bar{X}(t_1),\bar{X}(t_2),\cdots,\bar{X}(t_n))y(t_{i})\\
=
&\displaystyle \int \limits_{0}^{T}
\big{[}f_x(\bar{X}{(t)},\bar{u}(t))y(t)
+p(t)^{}(b(\bar{X}{(t)},u^{\varepsilon}(t))-b(\bar{X}{(t)},\bar{u}(t)))\big{]}dt.
\end{array}
\]
Therefore,
\begin{equation}%
\begin{array}
[c]{cl}
&
\displaystyle \sum \limits_{i=1}^{n}-\psi_{x_i}(\bar{X}(t_1),\bar{X}(t_2),\cdots,\bar{X}(t_n))y(t_{i})\\
= &
\displaystyle \int \limits_{0}^{T}
\big{[}f_x(\bar{X}{(t)},\bar{u}(t))y(t)
+p(t)^{}(b(\bar{X}{(t)},u^{\varepsilon}(t))-b(\bar{X}{(t)},\bar{u}(t)))\big{]}dt.
\end{array}
\label{maxm-2}%
\end{equation}
Now, let $u(t)=u$ be a constant, and note that $E_{\varepsilon}=[v,v+\varepsilon
]\subset \lbrack0,T]$. Combining equations (\ref{varm-3}) and (\ref{maxm-2}) and
noting the optimality of $\bar{u}(\cdot)$, we obtain%
\[%
\begin{array}
[c]{rl}%
0\leq & \hat{J}(u^{\varepsilon}(\cdot))-\hat{J}(\bar{u}(\cdot))\\
= & \displaystyle\sum_{i=1}^{n}\psi_{x_i}(\bar{X}(t_1),\bar{X}(t_2),\cdots,\bar{X}(t_n))y(t_i)\\
& +
{\displaystyle \int \limits_{0}^{T}}
\{f_x(\bar{X}{(t)},\bar{u}(t))y(t)+
f(\bar{X}{(t)},u^{\varepsilon}(t))-f(\bar{X}{(t)},\bar{u}(t))\}dt+o(\varepsilon)
\\
=&
-\displaystyle \int \limits_{0}^{T}
\{f_x(\bar{X}{(t)},\bar{u}(t))y(t)
+p(t)^{}(b(\bar{X}{(t)},u^{\varepsilon}(t))-b(\bar{X}{(t)},\bar{u}(t)))\}dt\\
&+{\displaystyle \int \limits_{0}^{T}}
\{f_x(\bar{X}{(t)},\bar{u}(t))y(t)+
f(\bar{X}{(t)},u^{\varepsilon}(t))-f(\bar{X}{(t)},\bar{u}(t))\}dt+o(\varepsilon)
\\
=&\displaystyle \int \limits_{0}^{T}
\{ H(\bar{X}(t),\bar{u}^{\varepsilon}(t),p(t))-H(\bar{X}(t),\bar{u}(t),p(t)) \}dt+o(\varepsilon).
\\
\end{array}
\]
Recalling
\begin{equation*}
H(x,u,p)=b(x,u)^{}p-f(x,u), \text{ \  \ }%
(x,u,p)\in \mathbb{R}^m\times U\times \mathbb{R}^m.%
\end{equation*}
Thus, we obtain
$$
H(\bar{X}(t),u,p(t))\leq H(\bar{X}(t),\bar{u}(t),p(t)),
$$
for any $u\in U$ and $t \in(t_{i-1},t_i)$, $i=1,2,\cdots,n$.

This completes the proof. $\ \ \ \ \ \ \ \ \Box$

\bigskip

In the following, we present the well-known Pontryagin's  maximum principle for the cost functional (\ref{cost-1}) under constrained conditions (\ref{ccc-1}). The cost functional is given as follows:
\begin{equation}
J(u(\cdot))=%
\big{[}{\displaystyle \int \limits_{0}^{T}}
f(X{(t)},u(t))dt+\Psi(X(T))\big{]},\label{cost-3}%
\end{equation}
the state process $X(\cdot)$ satisfies (\ref{ccc-1}).

\begin{theorem}
\label{ccc-th}
Let Assumptions (\ref{ass-b})-(\ref{ass-fai}) hold,
and let $(\bar{u}(\cdot),\bar{X}(\cdot))$ be an optimal pair of (\ref{cost-3}) under constrained conditions (\ref{ccc-1}).
Then there exists $(\beta^0,\beta^1,\cdots,\beta^n)\in \mathbb{R}^{n+1}$ satisfying
$$
\beta^0\geq 0,\ \ \left| \beta^0\right|^2+\displaystyle\sum_{j=1}^n\left| \beta^j\right|^2=1,
$$
and
$$
\displaystyle\beta^j(\gamma-\bar{X}(t_j))\geq 0,\  \gamma\leq 0,\ j=1,2,\cdots,n,
$$
and the adapted solution $p(\cdot)$ satisfying the following series of first-order adjoint equations,
\begin{equation}%
\begin{array}
[c]{ll}%
&-d{p}(t)= \{b_x(\bar{X}{(t)},\bar{u}(t))^{}p(t)-\beta^0f_x(\bar{X}{(t)},\bar{u}(t))\}dt,\ t\in(t_{i-1},t_{i}),\\
&p(t_{i})= -\beta^0\Psi_{x}(\bar{X}(t_n)1_{i=n}(i)-\beta^i+p(t_{i}^{+}),\ i=1,2,\ldots,n,
\end{array}
\label{prin-1}%
\end{equation}
and such that
\begin{equation}%
\begin{array}
[c]{ll}%
&H(\beta^0,\bar{X}(t),\bar{u}(t),p(t))\geq H(\beta^0,\bar{X}(t),u,p(t)),
\end{array}
\label{princ-3}%
\end{equation}
for any $u\in U$ and $t \in(t_{i},t_{i+1})$, $i=0,1,\cdots,n-1$, where
\begin{equation*}
H(\beta^0,x,u,p)=b(x,u)^{}p-\beta^0f(x,u),\text{ \  \ }
\end{equation*}
with $(\beta^0,x,u,p)\in \mathbb{R}\times\mathbb{R}^m\times U\times \mathbb{R}^m.$
\end{theorem}
\noindent\textbf{Proof}: Without loss of generality, we assume that $J(\bar{u}(\cdot))=0$, where $(\bar{u}(\cdot),\bar{X}(\cdot))$ is the optimal pair of problem (\ref{cost-3}) with constrained conditions (\ref{ccc-1}). For any $\theta>0$, we set
$$
J^{\theta}(u(\cdot))=\sqrt{\big{[}(J(u(\cdot))+\theta)^+\big{]}^2
+\displaystyle\sum_{i=1}^n\big{[}(-X^u(t_i))^+\big{]}^2}.
$$
From Assumption \ref{ass-fai}, one can verify that $J^{\theta}:\mathcal{U}[0,T]\to \mathbb{R}$ is continuous and satisfies
\begin{equation}
J^{\theta}(\bar{u}(\cdot))=\theta\leq \inf_{u\in\mathcal{U}[0,T]}J^{\theta}(u(\cdot))+\theta.
\end{equation}
Then, by Ekeland's variational principle, there exists a $u^{\theta}(\cdot)\in \mathcal{U}[0,T]$ such that
\begin{equation}
\label{ccci-ev}
J^{\theta}(u^{\theta}(\cdot))\leq J^{\theta}(\bar{u}(\cdot))=\theta,\ \tilde{d}(u^{\theta}(\cdot),\bar{u}(\cdot))\leq \sqrt{\theta},
\end{equation}
where $\tilde{d}(u^1(\cdot),u^2(\cdot))=M\{(t)\in[0,T]:u^1(t)\neq u^2(t)\}$, where $M$ is the product measure of the Lebesgue measure and probability on the set of $[0,T]$. We can also check that $(\mathcal{U}[0,T],\tilde{d})$ is a complete metric space. In addition, we have
\begin{equation*}
-\sqrt{\theta}\tilde{d}(u^{\theta}(\cdot),u(\cdot))\leq J^{\theta}(u(\cdot))-J^{\theta}(u^{\theta}(\cdot)),\ \forall u(\cdot) \in \mathcal{U}[0,T],
\end{equation*}
which deduces that
\begin{equation}
\label{cos-40}
J^{\theta}(u^{\theta}(\cdot))+\sqrt{\theta}\tilde{d}(u^{\theta}(\cdot),u^{\theta}(\cdot))\leq J^{\theta}(u(\cdot))+\sqrt{\theta}\tilde{d}(u^{\theta}(\cdot),u(\cdot)),\ \forall u(\cdot) \in \mathcal{U}[0,T].
\end{equation}
Thus, inequality (\ref{cos-40}) shows that $(u^{\theta}(\cdot),X^{\theta}(\cdot))$ is the optimal pair for the following cost functional
\begin{equation}
\label{cos-4}
J^{\theta}(u(\cdot))+\sqrt{\theta}\tilde{d}(u^{\theta}(\cdot),u(\cdot)),
\end{equation}
without the state constraints.

Because $U$ is a general control domain, let $\rho>0$ and $E_{\rho}=[v,v+\rho
]\subset (t_{i-1},t_i)$, for some $i\in\{1,2,\cdots,n\}$. Let $u\in U$ be any given constant. We define the following%
\[
u^{\theta,\rho}(t)=\left \{
\begin{array}
[c]{cl}%
{u}^{\theta}(t), & \text{if }t\in \lbrack0,T]\backslash E_{\rho},\\
u, & \text{if }t\in E_{\rho},
\end{array}
\right.
\]
which belongs to $ \mathcal{U}[0,T]$. It is easy to verify that
$$
\tilde{d}(u^{\theta,\rho}(\cdot),u^{\theta}(\cdot))\leq \rho.
$$
By equation (\ref{cos-40}), one obtains
\begin{equation}%
\begin{array}
[c]{rl}%
-\sqrt{\theta}\rho \leq & J^{\theta}(u^{\theta,\rho}(\cdot))-J^{\theta}(u^{\theta}(\cdot))\\
=&\displaystyle\frac{\big{[}(J(u^{\theta,\rho}(\cdot))+\theta)^+\big{]}^2-\big{[}(J(u^{\theta}(\cdot))+\theta)^+\big{]}^2
}
{J^{\theta}(u^{\theta,\rho}(\cdot))+J^{\theta}(u^{\theta}(\cdot))}\\
&+\displaystyle\frac{\sum_{j=1}^n\big{[}\big{[}(-X^{\theta,\rho}(t_j))^+\big{]}^2-
\big{[}(-X^{\theta}(t_j))^+\big{]}^2\big{]}}
{J^{\theta}(u^{\theta,\rho}(\cdot))+J^{\theta}(u^{\theta}(\cdot))},\\
\end{array}
\label{ccci-1}%
\end{equation}
where $X^{\theta,\rho}(\cdot))$ and $X^{\theta}(\cdot))$ are solutions of equation (\ref{ODE_1}) with controls $u^{\theta,\rho}(\cdot)$ and $u^{\theta}(\cdot)$, respectively. Let
\begin{equation}%
\begin{array}
[c]{ll}%
\beta^{0,\theta}=\displaystyle\frac{\big{[}J(u^{\theta}(\cdot))+\theta\big{]}^+}{J^{\theta}(u^{\theta}(\cdot))},\\
\beta^{j,\theta}=\displaystyle\frac{-\big{[}-X^{\theta}(t_j)\big{]}^+}
{J^{\theta}(u^{\theta}(\cdot))},\ j=1,2,\cdots,n.\\
\end{array}
\label{ccci-2}%
\end{equation}
Then, by the continuity of $J^{\theta}(\cdot)$ and Assumption \ref{ass-fai1}, we have
\begin{equation}%
\begin{array}
[c]{rl}%
& J^{\theta}(u^{\theta,\rho}(\cdot))-J^{\theta}(u^{\theta}(\cdot))\\
=&\beta^{0,\theta}\big{(}J(u^{\theta,\rho}(\cdot))-J(u^{\theta}(\cdot))\big{)}+
\displaystyle\sum_{j=1}^n\beta^{j,\theta}\big{(}X^{\theta,\rho}(t_j)
-X^{\theta}(t_j)\big{)}+o(\rho)\\
=&\displaystyle\sum_{j=1}^{n}\beta^{j,\theta}
(X^{\theta,\rho}(t_j)
-X^{\theta}(t_j))+\beta^{0,\theta}(\Psi(X^{\theta,\rho}(t_n))-
\Psi(X^{\theta}(t_n)))\\
&+\beta^{0,\theta}\displaystyle\int_0^T\big{(}f(X^{\theta,\rho}(t),u^{\theta,\rho}(t))-
f(X^{\theta}(t),u^{\theta}(t))\big{)}dt +o(\rho),
\end{array}
\label{ccci-3}%
\end{equation}
where $\frac{o(\rho)}{\rho}$ converges to $0$ when $\rho\to 0$.

Similar to Lemma \ref{lem-2}, let $(\bar{X}(\cdot),\bar{u}(\cdot))$ be replaced by $(X^{\theta}(t),u^{\theta}(t))$, $y(\cdot)$ be replaced by $\tilde{y}(\cdot)$ in equation (\ref{aprom-1}). Thus, one obtains
\begin{equation}%
\begin{array}
[c]{rl}
-\sqrt{\theta} \rho\leq & J^{\theta}(u^{\theta,\rho}(\cdot))-J^{\theta}(u^{\theta}(\cdot))\\
\leq& \beta^{0,\theta} \Psi_{x}(X^{\theta}(t_n))\tilde{y}(t_n) +
\displaystyle\sum_{i=1}^{n}\beta^{i,\theta}\tilde{y}(t_i)\\
&+\beta^{0,\theta}{\displaystyle \int \limits_{0}^{T}}
\big{\{}f_x(X^{\theta}{(t)},{u}^{\theta}(t))\tilde{y}(t)
+f(X^{\theta}{(t)},u^{\theta,\rho}(t))-f(X^{\theta}{(t)},{u}^{\theta}(t))\}dt+o(\rho).\\
\end{array}
\label{ccci-4}%
\end{equation}
In addition, we introduce the following adjoint equation,
\begin{equation}%
\begin{array}
[c]{rl}%
-d{p}^{\theta}(t)= & \{b_x(X^{\theta}(t),u^{\theta}(t))^{}p^{\theta}(t)
               -\beta^{0,\theta}f_x(X^{\theta}(t),u^{\theta}(t))\}dt,\ t\in(t_{i-1},t_{i}),\\
p^{\theta}(t_{i})= &-\beta^{0,\theta}\Psi_x(X^{\theta}(t_n))1_{i=n}(i)-
\displaystyle{\beta^{i,\theta}}{}+p(t_{i}^{+}), \text{ \ }i=1,\ldots,n.
\end{array}
\label{ccci-5}%
\end{equation}
Now, using the duality relation as in the proof of Theorem \ref{Maximumprinciple}, it follows that
\[%
\begin{array}
[c]{rl}%
 &o(1)+\sqrt{\theta} \geq E {\displaystyle \int \limits_{0}^{T}} \{H^{\theta,\rho}(t,{u}^{\theta,\rho}(t))-{H}^{\theta}(t,{X}^{\theta}{(t)})\}dt,\\
 \end{array}
\]
where
$$
{H}^{\theta}(t,{u}^{\theta}{(t)}):=H(\beta^{0,\theta},{X}^{\theta}{(t)},{u}^{\theta}(t),p^{\theta}(t))
$$
and
$$
H^{\theta,\rho}(t,{u}^{\theta,\rho}(t)):=H(\beta^{0,\theta},{X}^{\theta}{(t)},{u}^{\theta,\rho}(t),p^{\theta}(t)).
$$
Note that $o(1) \to 0$ when $\rho\to 0$. Thus, letting $\rho\to 0$, one obtains
\begin{equation}
\label{ccci-6}
\begin{array}
[c]{rl}%
 \sqrt{\theta}
 \geq &H^{\theta}(t,u(t))-{H}^{\theta}(t,u^{\theta}(t)).\\
 \end{array}
\end{equation}
From inequality (\ref{ccci-ev}), it follows that $u^{\theta}(\cdot)$ converges to $\bar{u}(\cdot)$ under distance $\tilde{d}$ as $\theta\to 0$. Then, by Assumptions \ref{ass-b}, \ref{assb-b2}, and the basic theory of differential equations, we have
$$
\displaystyle \sup_{0\leq t\leq T}\left|X^{\theta}(t)-\bar{X}(t)\right| \to 0,
$$
as $\theta\to 0$. By equation (\ref{ccci-2}), it follows that
\begin{equation}
\label{unq-1}
\left|\beta^{0,\theta}\right|^2+\displaystyle\sum_{j=1}^n\left|\beta^{j,\theta}\right|^2=1.
\end{equation}
Thus, we can choose a sequence $\{\theta_k\}_{k=1}^{\infty}$ satisfying $\displaystyle\lim_{k\to\infty}\theta_k=0$ and such that the limitations of $\beta^{0,\theta_k}$ and $\beta^{j,\theta_k}$ exist, and we use the notation
\begin{equation}
\begin{array}
[c]{ll}
\beta^{0}=\displaystyle\lim_{k\to\infty}\beta^{0,\theta_k},\\
\beta^{j}=\displaystyle\lim_{k\to\infty}\beta^{j,\theta_k},\\
\end{array}
\end{equation}
with $j=1,2,\cdots,n$. By equation (\ref{unq-1}), we have
$$
\left|\beta^{0}\right|^2+\displaystyle\sum_{j=1}^n\left|\beta^{j}\right|^2=1,
$$
and
$$
\displaystyle\beta^j(\gamma^j-\bar{X}(t_j))\geq 0,\  \gamma^j\geq 0,\ j=1,2,\cdots,n.
$$
Similarly, we can prove that
$$
\displaystyle \sup_{0\leq t\leq T}\big{[}\left|p^{\theta_k}(t)-p(t)\right|^2+\int_0^T\left|q^{\theta_k}(t)-q(t)\right|^2 \big{]} dt\to 0,
$$
as $k\to\infty$. Letting $k\to \infty$, from equation (\ref{ccci-6}), we have
\begin{equation}%
\begin{array}
[c]{ll}%
&H(\beta^0,\bar{X}(t),\bar{u}(t),p(t))\geq H(\beta^0,\bar{X}(t),u,p(t)),\\
\end{array}
\end{equation}
for any $u\in U$ and $t \in(t_{i},t_{i+1})$, $i=0,1,\cdots,n-1$.

Thus, we complete this proof. $\ \ \ \ \ \ \ \ \ \ \ \  \Box$

\bigskip

Similar to the proof of Theorem \ref{ccc-th}, we can obtain another kind of maximum principle for the multi-time state constraints.

\begin{corollary}\label{coro-1}
Let Assumptions (\ref{ass-b})-(\ref{ass-fai}) hold,
and let $(\bar{u}^n(\cdot),\bar{X}^n(\cdot))$ be an optimal pair of (\ref{cost-3}) under constrained conditions (\ref{ccc-1}).
Then, there exists $(\beta^{n,0},\beta^{n,1},\cdots,\beta^{n,n})\in \mathbb{R}^{n+1}$ satisfying
$$
\beta^{n,0}\geq 0,\ \ \left| \beta^{n,0}\right|^2+\displaystyle\sum_{i=1}^n\left| \frac{\beta^{n,i}}{n}\right|^2=1,
$$
and
$$
\displaystyle\beta^{n,i}(\gamma-\bar{X}^n(t_i))\geq 0, \  \gamma\leq 0,\ i=1,2,\cdots,n,
$$
and the adapted solution $p^n(\cdot)$ satisfying the following series of first-order adjoint equations,
\begin{equation}%
\begin{array}
[c]{ll}%
&-d{p}^n(t)= \{b_x(\bar{X}^n{(t)},\bar{u}^n(t))^{}p^n(t)
-\beta^{n,0}f_x(\bar{X}^n{(t)},\bar{u}^n(t))
-\displaystyle\frac{\sum_{j=i}^{n}\beta^{n,j}}{n}b_x(\bar{X}^n{(t)},\bar{u}^n(t))\}dt,\\
&p^n(t_{i})= -\beta^0\Psi_{x}(\bar{X}^n(t_n)1_{i=n}(i)+p^n(t_{i}^{+}),\text{ \ }t\in(t_{i-1},t_{i}),\ i=1,2,\ldots,n,
\end{array}
\label{4prin-1c}%
\end{equation}
and such that
\begin{equation}%
\begin{array}
[c]{ll}%
&H^n(\beta^{n,0},\bar{X}^n(t),\bar{u}^n(t),p^n(t))\geq H^n(\beta^{n,0},\bar{X}^n(t),u,p^n(t)),
\end{array}
\label{4princ-3c}%
\end{equation}
where
\begin{equation*}
H^n(\beta^{n,0},x,u,p)=b(x,u)^{}p-\beta^{n,0}f(x,u)
-\displaystyle\frac{\sum_{j=i}^{n}\beta^{n,j}}{n}b(x,u),\text{ \  \ }
\end{equation*}
for any $u\in U$ and $t \in(t_{i-1},t_{i})$, $i=1,2,\cdots,n$.
\end{corollary}
\textbf{Proof:} Similar to the proof of Theorem \ref{ccc-th}, for any $\theta>0$, we use the following cost functional without constrained conditions,
$$
J^{n,\theta}(u(\cdot)):=\sqrt{\big{[}(J(u(\cdot))+\theta)^+\big{]}^2
+\displaystyle\sum_{i=1}^n\frac{1}{n^2}\big{[}(-X^u(t_i))^+\big{]}^2}.
$$
By Ekeland’s variational principle, we can show that there exists $(u^{n,\theta}(\cdot),X^{n,\theta}(\cdot))$, which is the optimal pair of the following cost functional
\begin{equation}
\label{cos-4c}
J^{n,\theta}(u(\cdot))+\sqrt{\theta}\tilde{d}(u^{n,\theta}(\cdot),u(\cdot)),
\end{equation}
without the state constraints, where $\tilde{d}(u^1(\cdot),u^2(\cdot))=M\{(t)\in[0,T]:u^1(t)\neq u^2(t)\}$.

We define the following%
\[
u^{n,\theta,\rho}(t)=\left \{
\begin{array}
[c]{cl}%
{u}^{n,\theta}(t), & \text{if }t\in \lbrack0,T]\backslash\ E_{\rho},\\
u, & \text{if }t\in E_{\rho},
\end{array}
\right.
\]
which belongs to $ \mathcal{U}[0,T]$. Let
\begin{equation}%
\begin{array}
[c]{ll}%
\beta^{n,0,\theta}=\displaystyle\frac{\big{[}J(u^{n,\theta}(\cdot))+\theta\big{]}^+}{J^{n,\theta}(u^{n,\theta}(\cdot))},\\
\beta^{n,i,\theta}=\displaystyle\frac{-\big{[}-X^{n,\theta}(t_i)\big{]}^+}
{J^{n,\theta}(u^{n,\theta}(\cdot))},\ i=1,2,\cdots,n,\\
\end{array}
\label{ccci-2c}%
\end{equation}
where $X^{n,\theta,\rho}(\cdot))$ and $X^{n,\theta}(\cdot))$ are the related solutions of equation (\ref{ODE_1}) with controls $u^{n,\theta,\rho}(\cdot)$ and $u^{n,\theta}(\cdot)$. Then, by the continuity of $J^{n,\theta}(\cdot)$ and Assumption \ref{ass-fai}, we have
\begin{equation}%
\begin{array}
[c]{rl}%
& J^{n,\theta}(u^{n,\theta,\rho}(\cdot))-J^{n,\theta}(u^{n,\theta}(\cdot))\\
=&\beta^{n,0,\theta}\big{[}J(u^{n,\theta,\rho}(\cdot))-J(u^{n,\theta}(\cdot))\big{]}+
\displaystyle\sum_{i=1}^n\beta^{n,i,\theta}\big{[}X^{n,\theta,\rho}(t_i)
-X^{n,\theta}(t_i)\big{]}+o(1),\\
=&\displaystyle\sum_{i=1}^{n}\frac{\beta^{n,i,\theta}}{n}
\left[X^{n,\theta,\rho}(t_i)
-X^{n,\theta}(t_i)\right]+\beta^{n,0,\theta}\left[\Psi(X^{n,\theta,\rho}(t_n))-
\Psi(X^{n,\theta}(t_n))\right]\\
&+\beta^{n,0,\theta}\displaystyle\int_0^T\big{[}f(X^{n,\theta,\rho}(t),u^{n,\theta,\rho}(t))-
f(X^{n,\theta}(t),u^{n,\theta}(t))\big{]}dt  +o(\rho).
\end{array}
\label{ccci-3c}%
\end{equation}
Note that
$$
X^{n,\theta,\rho}(s)
-X^{n,\theta}(s)=\int_0^s\left[b(X^{n,\theta,\rho}(t),u^{n,\theta,\rho}(t))-b(X^{n,\theta}(t),u^{n,\theta}(t))\right]dt,
$$
it follows that
\begin{equation}%
\begin{array}
[c]{rl}%
  & \displaystyle\sum_{i=1}^{n}\frac{\beta^{n,i,\theta}}{n}
\left[X^{n,\theta,\rho}(t_i)
-X^{n,\theta}(t_i)\right]\\
=&\displaystyle\sum_{i=1}^{n}\int_{t_{i-1}}^{t_i}\frac{\sum_{j=1}^{n-i+1}\beta^{n,j,\theta}}{n}
\left[b(X^{n,\theta,\rho}(t),u^{n,\theta,\rho}(t))-b(X^{n,\theta}(t),u^{n,\theta}(t))\right]dt.
\end{array}
\label{ccci-3c}%
\end{equation}
The last step is the same as in the proof of Theorem \ref{ccc-th}; therefore, we omit it.  $\ \ \ \ \ \ \ \ \Box$

\end{document}